\input amstex
\mag=\magstep1
\documentstyle{amsppt}
\topmatter
\title On Classification of the Extremal Contraction from a Smooth Fourfold
\endtitle
\author Hiromichi Takagi\endauthor
\rightheadtext {extremal contraction}
\address Department of Mathematical Science,Tokyo University,Komaba,Tokyo,153,Japan\endaddress
\abstract We classify extermal divisorial contraction which contracts a divisor
to a curve from a smooth fourfold.We prove the exceptional divisor is 
$\Bbb P^2$bundle or quadric bundle over a smooth curve and the contraction is the blowing up along the curve.\endabstract
\endtopmatter

\document

\head 0.Introduction\endhead
In his pioneer paper \cite{M1} and \cite{M2}, Shigefumi Mori introduced 
the extremal ray and classified completely the extremal contraction 
from a smooth 3-fold. In dimension 4, we will consider the same problem, i.e., we want to classify the extermal contraction from a smooth 4-fold. In dimension 4, the situation is more complicated.

First small (flipping) contraction appears. This case was completely classified 
by Yujiro Kawamata in his ingeneous paper\cite{Ka1}.

Secondly, in case the contraction is fibre type from 4-fold to 3-fold and 
divisorial type which contracts a divisor to a surface, equidimentionality
of the fibre is not satisfied in general.(i.e., general fibres are 1 dimensional
but some special fibres are possibly 2 dimensional.) The special 2 dimensional 
fibre are classified by Yasuyuki Kachi in \cite{Kac} in case of fibre type, and
by Marco Andreatta in case of divisorial type.

Thirdly in case the contraction is divisorial type which contracts a divisor to a point, the exceptional divisor is possibly nonnormal. In fact Mauro Beltrametti listed up all the possibility of the exceptional divisor in \cite{Be1},
\cite{Be2}. It include nonnormal possibility.(But many cases are excluded by 
Takao Fujita in \cite{F2}.)

In this paper, we consider the contraction is divisorial type which contracts 
a divisor to a curve. This case turns out to be very mild in contrast to the above cases. 

\proclaim{Main Theorem}
Let $X$ be a smooth 4-fold and let $f\:X\to Y$ be a divisorial contraction which contracts a divisor to a curve. 
Let $E$ be the exceptional divisor of $f$ and $C$ be $f(E)$.
Then
\roster
\item"(1)" $C$ is a smooth curve .
\item"(2)" $f\vert_E\:E\to C$ is $\Bbb P^2$bundle or quadric bundle(see the difinition1.4 below) over $C$.
\item"(3)" $f$ is the blowing up of $Y$ along $C$.
\endroster
\endproclaim
\definition{Acknowledgement}
I express my gratitude to Professor Y.Kawamata for encouraging me during
preparing this paper. In particular, he told me that I should prove the irreducibility of general fibres in Theorem1.1 and investigate the local structure of the contraction. I am also thankful to Dr.Yoshiaki Fukuma for giving me much advice.
\enddefinition

\head 1.Notations and Preliminaries\endhead 
We cite the key theorems and make some definitions in this section.
\definition{Notation1.0}The $\Bbb P^1$ bundle $\Bbb P(\Cal O_{\Bbb P^1}\oplus \Cal O_{\Bbb P^1}(-a))$ over $\Bbb P^1$ is called a Hirzebruch surface of degree $a$ and denoted by $\Bbb F_a$. The unique negative section of it is denoted by $C_0$ and a ruling is denoted by $f$. 

The projective cone obtained from $\Bbb F_a$ by the contraction of $C_0$ is denoted by $\Bbb F_{a,0}$. A generating line on $\Bbb F_{a,0}$ is denoted by $l$.\hfill $\square$ \enddefinition  
\proclaim{Theorem1.1}(cf.\cite{TA},\cite{Be1} and \cite{Be2})
Let $X$,$Y$,$E$ and $C$ be as in Main Theorem.Let $F$ be a general fibre of $f\vert_E\:E\to C$.

Then 
$$(F,-K_X\vert_F) \simeq(\Bbb P^2,\Cal O_{\Bbb P^2}(1)),(\Bbb P^2,\Cal O_{\Bbb P^2}(2)),
(\Bbb P^1 \times \Bbb P^1,\Cal O(1,1)) or (\Bbb F_{2,0}, \Cal O_{\Bbb P^3}(1)\vert_{\Bbb F_{2,0}})$$
\hfill $\square$
\endproclaim

We give the outline of the proof.
\demo{Outline of the proof} 
Once we prove the irreducibility of $F$, the results are follow from [TA] and [Be1] and [Be2], so we will prove only the irreducibility here.     
We assume that $F$ is reducible and get a contradiction.
Let $H$ be a good supporting divisor of $f$. We may assume that $H$ is a smooth variety and at least locally $F=H\cap E$. By the adjunction formula, $-K_F=-K_X-E-H\vert_F$, but since $H\vert_F\sim0$, $-K_F=-K_X-E\vert_F$. Note that $-K_X\vert_F$ and $-E\vert_F$ is ample Cartier divisor on $F$. Since $-K_F$ is ample on $F$, $F$ is a generalized del Pezzo surface(i.e., a Gorenstein (possibly reducible) anti polarized surface). Write $F=\cup F_i$ where $F_i$ is a irreducible component of $F$. By [R3], $(F_i,-K_F\vert_{F_i})$ is one of the following:
\roster
\item"(a)" $(\Bbb P^2, \Cal O_{\Bbb P^2}(i))$ where $i$ is 1 or 2.       
\item"(b)" $(\Bbb F_{a,0}, \Cal O_{\Bbb F_{a,0}}(al))$
\item"(c)" $(\Bbb F_a, \Cal O_{\Bbb F_a}(C_0+(a+1)f))$
\item"(d)" $(\Bbb F_a, \Cal O_{\Bbb F_a}(C_0+(a+2)f))$
\endroster But since $-K_F$ is a sum of two ample Cartier divisor, (a) with $i=1$,(b),(c) and (d) is impossible. So we get $(F_i, -K_F\vert_{F_i})\simeq(\Bbb P^2, \Cal O_{\Bbb P^2}(2))$ and $-K_X\vert_{F_i}\simeq \Cal O_{\Bbb P^2}(1)$. Furthermore by[R3,Main Theorem and 1.3], $F$ is union of two $\Bbb P^2$'s which intersect line in $\Bbb P^2$. But this is impossible since $N_{{F_i}/ H}\simeq \Cal O_{\Bbb P^2}(-2)$ and so the birational contraction which contracts only one $F_i$.\hfill $\square$ \enddemo

The next theorem of freeness is very useful for classification of low 
dimensional fibres of an extremal contraction.
\proclaim{Theorem1.2}(see \cite{A-W}) Let $X$ be a normal log terminal variety
and $L$ be an ample line bundle on $X$. Let $f\:X\to Y$ be the adjoint 
contraction supported by $K_X+rL$ and $F$ be a fibre of $f$.
Assume that $\dim F<r+1$ if $\dim Y<\dim X$ or $\dim F\leq r+1$ if
$\dim Y\leq \dim X$.

Then $f^*f_*L\to L$ is surjective at every point of $F$.
\hfill $\square$
\endproclaim
The next theorem was proved in \cite{Wi1} and \cite{Wi2}(see also \cite{R1} and \cite{R2}) but we will give the proof again here for readers'convinience.
\proclaim{Theorem1.3}Let $X$ be a smooth 3-fold and $Y$ be a canonical 3fold.
Let $f\:X\to Y$ be a crepant birational contraction which contracts a 
irreducible divisor to a curve. Let $E$ be the exceptional divisor and $C$ be
$f(E)$. 

Then $C$ is a smooth curve.
\endproclaim
\demo{Proof}Let $P$ be any point of $C$.
The assertation is local, so we may replace $C$ and $Y$ with an affine(not analytic) neighborhood of $P$. We will keep this in mind below.
\proclaim{Claim1}$P$ is a $cDV$ point of $Y$. 
\endproclaim

\demo{Proof}Suppose that $P$ is not a $cDV$ point. By \cite{R1} and \cite{R2}, we have a birational morphism $g\:X'\to X$ such that $X'$ is terminal, $g$ is crepant and $g$ has a exceptional divisor $E_0$ which contracts to $P$. Take a common resolution $\tilde X$.$$\CD 
{\tilde X} @> >> X' \\
@VVV @VVgV \\
X @>> f> Y
\endCD$$ 
Since $f$ and $g$ is crepant, strict transform of $E_0$ on $X$ is exceptional for $f$ but this contradicts to the irreducibility of the exceptional divisor of $f$.\hfill $\square$ \enddemo
 Let $H$ be the pull back of a very ample divisor on $Y$.
\proclaim{Claim2}$\vert mH-E\vert$ is very ample for $m\gg 0$. In particular $f^*f_*\Cal O_X(-E)\to \Cal O_X(-E)$ is surjective.\endproclaim
\demo{Proof}
Since $mH-2E$ is ample for $m\gg 0$, it follow from the vanishing theorem
(see \cite{KMM}) and the exact sequence $$ 0\to \Cal O_X(mH-2E)\to \Cal O_X(mH-
E)\to \Cal O_E(mH-E)\to 0$$ that $H^0(\Cal O_X(mH-E))\to H^0(\Cal O_E(mH-E))$ is surjective. Let $l$ be a fibre of $f\vert_E$. From the vanishing $H^1(X,\Cal O_X)=0$, $l$ is a tree of $\Bbb P^1$, so $\vert mH-E\vert_l\vert$ is very ample and so is $\vert mH-E\vert_E\vert$ since we consider locally. From thease, $\vert mH-E\vert$ is also very ample.\hfill $\square$ \enddemo

\proclaim{Claim 3}$C$ can be embedded in a smooth surface.\endproclaim
\demo{Proof}
In fact, let $S$ be a smooth general member of $\vert mH-E \vert$. Since $f\vert_S\: S\to f(S)$ is e\'tale, so $f(S)$ is smooth. $C$ is in $f(S)$ so we are done.\hfill$\square$ \enddemo  
\proclaim{Claim4}$f$ is the blowing up of $Y$ along $C$.
\endproclaim
\demo{Proof}Let $\Cal I_C$ be the ideal sheaf of $C$ in $Y$. 

First we see that$$f_*\Cal O_X(-E)=\Cal I_C$$ Let's consider the exact sequence$$0\to \Cal O_X(-E)\to \Cal O_X\to \Cal O_E\to 0$$ From this and vanishing theorem, we get the exact sequence$$0\to f_*\Cal O_X(-E)\to \Cal O_Y\to f_*\Cal O_E\to 0$$ Since except at finite points $f_*\Cal O_E=\Cal O_C$, we have $f_*\Cal O_X(-E)=\Cal I_C$ except at finite point. But they are reflexive(cf.[H2,Corollary 1.5,Proposition 1.6,Corollary 1.7]),so they are actually equal.     

>From this and claim2, we have $f^*\Cal I_C\to \Cal O_X(-E)$ is surjective, i.e., $\Cal I_C\Cal O_X=\Cal O_X(-E)$. So by the universal property of blowing up(cf.[H1,II proposition7.14]), $f$ decomposes as $X\to X_1\to Y$, where $X_1$ is the blowing up of $Y$ along $C$. $X_1$ is normal since $C$ is a $cDV$ curve in $Y$.(see also the calculations below.) So $X\simeq X_1$ since $f$ is a primitive contraction. This established claim4.
\hfill $\square$
\enddemo

We suppose $P$ is a singular point of $C$ and get a contradiction.
Recall that $P$ is a $cDV$ point and so $Y$ can be embedded(analytically locally)in $\Bbb C^4$. Let $x,y,z,t$ be its coordinate around $P$ and $g$ be the defining equation of $Y$ in $\Bbb C^4$. By claim3, We may assume $\Cal I_C=<x,y,h>$, where $h\in \Bbb C[[z,t]]$. Since we suppose $P$ is singular point of $C$, $h=0$ is singular at $z=t=0$ in $zt$-plane. $X$ is the strict transform of $Y$ in the blowing up $\tilde \Bbb C^4$ of $\Bbb C^4$ along $C$. $\tilde \Bbb C^4$ is in $\Bbb C^4\times \Bbb P^2$ and given by $$\operatorname{rank}\pmatrix
x & y & h \\
u & v & w   
\endpmatrix \leq 1$$ ($u,v,w$ is the homogenious coordinate of $\Bbb P^2$).
Take the affine piece given by $u=1$. We can embed this affine piece of $\tilde \Bbb C^4$ in $\Bbb C^5$ with coordinate $(x,z,t,v,w)$ and equation $xw=h$. This affine variety is singular above $P$ along the line $L$ defined by $x=z=t=w=0$. Write $g(x,xv,z,t)=x^m \tilde g(x,v,z,t)$,where $\tilde g$ can not devided by $x$. Then $X$ is defined by $xw=h$ and $\tilde g=0$. This intersects $L$ and at the intersections, $X$ is singular, a contradiction. 
\hfill $\square$    
\enddemo
\definition{Definition1.4}Let $E$ be a normal projective 3-fold and $C$ be a smooth curve. Let $f\: E\to C$ be a projective surjective morphism.

We say $f\: E\to C$ is quadric bundle if the following conditions are satisfied.\roster
\item"(1)" there exists a $f$-very ample line bundle $\Cal L$ on $E$. 
\item"(2)" For any closed point $s$ of $C$, $h^0(E_s,\Cal L_s)=3$ and $E_s$ is a quadric in $$\Bbb P(H^0(E_s,\Cal L_s)^*)\simeq \Bbb P^3$$ 
\endroster
If we can take such an $\Cal L$, we say $f\: E\to C$ is the quadric bundle associoated to $\Cal L$.    
\hfill $\square$
\enddefinition
\definition{Remark}We require $E$ is normal above so general fibres of $f$ are irreducible.\hfill $\square$ \enddefinition  
\head 2.Proof of the Main theorem\endhead

\demo{Proof of (1)}Let $H$ be a good supporting divisor of $f$ and $L$ be $\Cal O_X(mH-K_X)$ for $m\gg0$. Since $L$ is ample, we apply Theorem1.2 for this $f$ and $L$ with $r=1$.(Remark that the dimension of fibres of $f$ is 2.) Then for any point $P$ of $C$ and a suitable affine neighborhood $U$ of $P$ in $Y$, $\vert L\vert_{f^{-1}(U)}\vert$ has no base points on $E$. So furthermore if we replace $f^{-1}(U)$ with a suitable neighborhood $V$ of $E$ in $X$, $\vert L\vert_V\vert$ is base point free on $V$. Let $X_0$ be a general smooth member of $\vert L\vert_V\vert$, $Y_0$ be $f(X_0)$ and $E_0$ be $E\vert_{X_0}$. Then     
\proclaim{Claim}$E_0$ is irreducible \endproclaim
\demo{Proof}Let $F$ be a general fibre of $f$ and $F_0$ be $F\vert_{X_0}$ . It suffices to prove $F_0$, i.e., $X_0\vert_F$  is irreducible. For this it suffices to prove the surjectivity of $H^0(L\vert_V)\to H^0(L\vert_F)$. For by Theorem1.1, general member of $\vert L\vert_F\vert$ is irreducible. First from the exact sequence $$0\to L\vert_V\otimes \Cal O_X(-E)\to L\vert_V\to L\vert_E\to 0$$ and the vanishing theorem, we have $H^0(L\vert_V)\to H^0(L\vert_E)$ is surjective. Secondly since we take $F$ to be a general fibre, $F$ is a Cartier of $E$ and since near $F$, $E\simeq \Bbb P^2\times \Bbb A^1, \Bbb P^1\times \Bbb P^1\times \Bbb A^1 or \Bbb F_{2,0}\times \Bbb A^1$, $E$ has at worst canonical singularities near $F$(cf.Theorem1.1). So we can use the vanishing theorem for the exact sequence $$0\to L\vert_E \otimes \Cal O_E(-F)\to L\vert_E\to L\vert_F\to 0$$ and we get the surjectivity of $H^0(L\vert_E)\to H^0(L\vert_F)$. This establishes the claim.\hfill $\square$ \enddemo 

>From this claim, we can apply Theorem1.3 for $X_0$, $Y_0$, $E_0$ and prove that $C$ is smooth.\hfill $\square$ \enddemo

\demo{Proof of (2)}Let $F$ be a general fibre of $f$. If $(F, -K_X\vert_F)\simeq (\Bbb P^2, \Cal O_{\Bbb P^2}(1)), (\Bbb P^1\times \Bbb P^1, \Cal O(1,1)), or (\Bbb F_{2,0}, \Cal O_{\Bbb P^3}(1)\vert_{\Bbb F_{2,0}})$, let $\Cal L$ be $\Cal O_E(-K_X)$. 

If $(F, -K_X\vert_F)\simeq (\Bbb P^2, \Cal O_{\Bbb P^2}(2))$, let $\Cal L$ be $\Cal O_E(-E)$. Then we will prove \roster
\item"({\bf i})" If $F\simeq \Bbb P^2$, $f\vert_E \: E\to C$ is $\Bbb P^2$bundle.
\item"({\bf ii})" If $F\simeq \Bbb P^1\times \Bbb P^1 or \Bbb F_{2,0}$, $f\vert_E\: E\to C$ is the quadric bundle associated to $\Cal L$.
\endroster 
If $\Cal L \simeq \Cal O_E(-K_X)$, we can argue as follow.

First we see $f_*\Cal L$ is locally free. The exact sequence$$0\to \Cal O_X(-K_X-E)\to \Cal O_X(-K_X)\to \Cal O_E(-K_X)\to 0$$ and the vanishing theorem, we have $R^i{f_*\Cal O_E(-K_X)}=0$($i>0$). On the otherhand, we have $H^i(E_s, \Cal L_s)=0$ for $i\gg0$ and any $s\in C$. Furthermore by above (1), $C$ is smooth so $f\vert_E$ is flat. So by Cohomology and Base change theorem(cf.[H1,III Theorem12.11]), $f_*\Cal L$ is locally free.     

Next we see $\Cal L$ is $f$-free. Let's consider the commutative diagram$$\CD
f^*f_*\Cal O_X(-K_X) @> >> f^*f_*\Cal O_E(-K_X)\\ 
@VVV @VVV \\
\Cal O_X(-K_X) @> >> \Cal O_E(-K_X)\endCD$$
We see the left arrow is surjective by Theorem1.2 and so is the bottom arrow by the above exact sequence and vanishing. So the right arrow must be surjective, i.e., $\Cal L$ is $f$-free.

By thease, we get the morphism $g\: E\to \Bbb P (f_*\Cal L)$ defined by $\Cal L$.

If we are in case(i), $g$ is birational since general fibre of $f$ is $\Bbb P^2$ and $\Bbb P(f_*\Cal L)$ is $\Bbb P^2$bundle. And $g$ is finite since $\Cal L$ is $f$-ample. So by Zariski Main Theorem, $g$ is isomorphism.

If we are in case(ii), then 
\proclaim{Claim}$\Cal L$ is $f$-very ample. 
\endproclaim \demo{Proof}We will prove $g$ is isomorphism onto $g(E)$. $g$ is birational because on the general fibre of $f$, $\Cal L$ is very ample. $g$ is finite because $\Cal L$ is $f$-ample. The dimension of singular locus of $g(E)$ is not greater than 2 since general fibres of $g(E)\to C$ are $\Bbb P^1\times \Bbb P^1$ or $\Bbb F_{2,0}$ and $C$ is smooth. $g(E)$ satisfies the Serre's condition $S_2$ since $E$ is a divisor of smooth 4-fold. So $g(E)$ is normal. Then by Zariski main theorem, $g$ is isomorphism.\hfill $\square$ \enddemo From this claim, it is easy to see that $E$ is quadric bundle associated to $\Cal L$.

If $\Cal L\simeq \Cal O_E(-E)$, we can argue as follow.(cf.[F1,1.5])

Let $F'$ be any fibre of $f\vert_E$ and write $F'=\cup F'_i$, where $F'_i$ is a irreducible component of $F'$. Since $f\vert_E$ is flat, $1=(-E)^2 F=\sum(-E)^2F_i$. Since $-E\vert_{F_i}$ is ample, $(-E)^2 F_i>0$. So $F'$ must be irreducible. By the lower semicontinuity of $\Delta$-genus(cf.[H1,III,Theorem 12]), we have $\Delta (F',-E\vert_{F'})\leqq \Delta (F,-E\vert_F)=0$. $F'$ has no embedded points since $E$ is Cohen-Macauley and $F'$ is Cartier divisor on $E$. So $\Delta(F',-E\vert_{F'})\geqq0$ by [F0] and so $\Delta(F',-E\vert_{F'})=0$. Since $(-E)^2 F'=1$, $(F',-E\vert_{F'})\simeq (\Bbb P^2,\Cal O(1))$ by the classification of the varieties of $\Delta$-genus 0. So $f\: E\to C$ is $\Bbb P^2$bundle.     
\hfill $\square$
\enddemo
\definition{Remark}We cannot proceed in case $\Cal L\simeq \Cal O_E(-E)$ similar to the case $\Cal L\simeq \Cal O_E(-K_X)$ because we have not freeness of $\Cal O_E(-E)$ apriori.\hfill $\square$ \enddefinition    

As for (3), the proof is almost the same as [M2,Corollary(3.4)]. So we will give only the outline of the proof.
\demo{Outline of the proof of (3)}We see that $\Cal O_E(-E)$ is $f\vert_E$-very ample by (2) and $\Cal O_X(-E)$ is $f$-ample. So we get the following.
\proclaim{Claim}
\roster
\item"(a)" $R^i f_* \Cal O_X(-jE)=0$ for $i>0$ and $j\geq0$.  
\item"(b)" $f_* \Cal O_X(-jE)={\Cal I_C}^j$, ${\Cal I_C}^j\Cal O_X=\Cal O_X(-jE)$ for $j\geq0$.
\item"(c)" $\oplus_{n\geq 0} {\Cal I_C}^n /{\Cal I_C}^{n+1}\simeq \oplus_{n\geq 0}f_* \Cal O_E(-nE)$ as $\Cal O_C$ algebra.
\endroster \hfill $\square$ \endproclaim 
By this claim, we can easily get the result.\hfill $\square$ \enddemo 

\definition{Remarks and Examples}We can say the following about the local analytic structure of the contraction. Let $F'$ be a fibre of $f\vert_E$ and $F$ is any general fibre of $f\vert_E$ near $F'$. We will give the description near $F'$.
\roster
\item"(1)" If $f\vert_E\: E\to C$ is $\Bbb P^2$ bundle and $\Cal O_X(-E)\vert_F\simeq \Cal O_{\Bbb P^2}(1)$, $Y$ is smooth along $C$.(cf.[SN]) 
\item"(2)" If $f\vert_E\: E\to C$ is $\Bbb P^2$ bundle and $\Cal O_X(-E)\vert_F\simeq \Cal O_{\Bbb P^2}(2)$, $Y$ can be considered as one parameter family of $\frac12(1,1,1)$ singularity. In fact, let $P$ be any point of $C$ and take a general very ample divisor $A$ through $P$. Let $H$ be the pull back of $A$. Then $H$  is smooth along $E\vert_H\simeq \Bbb P^2$ since $E\vert_H$ is smooth and a Cartier divisor of $H$. So $f\vert_H$ is the extremal contraction from a smooth 3-fold near the fibre over $P$. Then by Mori's classification, $A$ has $\frac12(1,1,1)$ singularity at $P$.      
\item"(3)" If $f\vert_E\: E\to C$ is quadric bundle we see $Y$ is locally hypersurface in $\Bbb C^5$ and $C$ is locally complete intersection in the $\Bbb C^5$
 because $\Cal I_C/{\Cal I_C}^2$ is locally free sheaf of rank 4 on $C$ by the claim in the proof of Main Theorem (3). But the type of singularity of $Y$ along $C$ is very various(case (c) and (e) below), so we will give some examples here. Below $\Bbb C^5$ has always coordinates $x,y,z,w,t$. $Y$ is given hypersurface in $\Bbb C^5$ and $X$ is the blow up of $Y$ along $C$. We assume $F'$ is the fibre over the origin. 
\item"(3a)"($F$ and $F'$ are $\Bbb P^1\times \Bbb P^1$) Let $Y$ be $(x^2+y^2+z^2+w^2=0)$ and $C$ be $(x=y=z=w=0)$. In this case above example is all.    
\item"(3b)"($F$ is $\Bbb P^1\times \Bbb P^1$ and $F'$ is $\Bbb F_{2,0}$) Let $Y$ be $(x^2+y^2+z^2+tw^2=0)$ or $(x^2+y^2+z^2+t^m w^2+w^3=0)$ and $C$ be $(x=y=z=w=0)$. In this case above examples are all.
\item"(3c)"($F$ is $\Bbb P^1\times \Bbb P^1$ and $F'$ is union of two planes in $\Bbb P^3$) Let $Y$ be $(x^2+y^2+t^m z^2+z^3+t^n w^2+w^3=0)$ and $(x=y=z=w=0)$. 
\item"(3d)"($F$ and $F'$ are $\Bbb F_{2,0}$)  Let $Y$ be $(x^2+y^2+z^2+w^3=0)$ and $(x=y=z=w=0)$. In this case this example is all.
\item"(3e)"($F$ is $\Bbb F_{2,0}$ and $F'$ is union of two planes in $\Bbb P^3$)  Let $Y$ be $(x^2+y^2+z^3+t^m w^2+w^3=0)$ and $(x=y=z=w=0)$. 
\endroster \hfill $\square$ \enddefinition

\definition{Question} Dose any quadric bundle appear as the exceptional divisor of the contraction as in Main Theorem ?
\enddefinition

\enddocument